
\magnification=1200
\hfuzz=1pt
\def\sqr#1#2{{\vcenter{\vbox{\hrule height.#2pt \hbox{\vrule width.#2pt height#1pt \kern#1pt \vrule width.#2pt} \hrule height.#2pt}}}}
\def\square{\mathchoice\sqr64\sqr64\sqr{2.1}3\sqr{1.5}3}
\def\forc{{\,\parallel\joinrel\relbar\joinrel\relbar\,}}
\def\qed{\hfill{\vbox{\hrule\hbox{\vrule\kern3pt
                \vbox{\kern6pt}\kern3pt\vrule}\hrule}}}
\def\Dom{{\rm Dom}}
\def\Ran{{\rm Ran}}
\def\cf{{\rm cf}}
\def\supp{{\rm supp}}

\centerline{\bf Universal graphs without  large cliques}
\medskip  
\centerline{\bf  P. Komj\'ath, S.
Shelah\footnote{$^*$}{{\rm Publication No. 492. Research partially
supported by BSF.}}}
\medskip
\centerline{\sl Dedicated to the memory of Alan Mekler}
\bigskip
\centerline{\bf 0. Introduction}
\medskip
\noindent The theory of universal graphs originated 
from the observation of R.~Rado
[4,5] that a universal countable graph $X$ exists, i.e., $X$ is countable and
isomorphically embeds every countable graph. He also showed that under GCH,
there is a universal graph in every infinite cardinal. Since then, several
results have been proved about the existence of universal elements 
 in different classes of graphs. 
For example, a construction similar to Rado's shows, that for every natural
number $n\geq 3$, there is a universal $K(n)$-free countable graph, or, if
GCH is assumed, there is one in every infinite cardinal (here $K(n)$ 
denotes the
complete graph on $n$ vertices). This result also
follows from the existence theorem of universal and special models.
 
The following folklore observation shows that this cannot be extended to
$K(\omega)$. Assume that $X=(V,E)$ is a $K(\omega)$-free graph of cardinal
$\lambda$ that embeds every $K(\omega)$-free graph of cardinal $\lambda$.
Let $a\not\in V$, and define the graph $X'$ on $V'=V\cup\{a\}$ as follows.
$X'$ on $V$ is identical with $X$,  $a$ is joined to every vertex
of $V$. Clearly, $X'$ is $K(\omega)$-free. So, by assumption, there is an
embedding $g\colon V'\to V$ of $X'$ into $X$. Put $a_0=a$, and, by
induction, $a_{n+1}=g(a_n)$. As $g$ is edge preserving, we get, by
induction on $n$, that $a_n$ is joined to every $a_t$ with $t>n$, so they
are distinct, and form a $K(\omega)$ in $X'$, a contradiction.
 
In Section 1 we give some existence/nonexistence statements on universal
graphs, which under GCH give a necessary and sufficient condition for the
existence of a universal graph of size $\lambda$ with no $K(\kappa)$,
namely, if either $\kappa$ is finite or $\cf (\kappa) > \cf (\lambda)$. The
special case when $\lambda^{<\kappa}=\lambda$ was first proved by
F.~Galvin.
 
In Section 2 we investigate the question that if there is no
universal $K(\kappa)$-free graph of size $\lambda$ then how many of these
graphs embed all the other. It was proved in [1], that if
$\lambda^{<\lambda}=\lambda$ (e.g., if $\lambda$ is regular and the GCH
holds below $\lambda$), and $\kappa=\omega$, then this number is
$\lambda^+$. We show that this holds for every $\kappa\leq\lambda$ of
countable cofinality. On the other hand, even for $\kappa=\omega_1$, and
any regular $\lambda\geq\omega_1$ it is consistent that the GCH holds below
$\lambda$, $2^{\lambda}$ is as large as we wish, and the above number is either
$\lambda^+$ or $2^{\lambda}$, so both extremes can actually
occur. Similar results when the excluded graphs are
disconnected, were proved in [2] and [3].
 
\medskip
\noindent{\bf Notation.} We use the standard axiomatic set
theory notation. If $X$ is a set, $\kappa$ a cardinal,
$[X]^{\kappa}=\{Y\subseteq X\colon \vert Y\vert =\kappa\}$, 
$[X]^{<\kappa}=\{ Y\subseteq X\colon \vert Y\vert < \kappa\}$. A {\sl
graph} is a pair $X=(V,E)$ where $V$ is some set, and $E\subseteq [V]^2$,
i.e., we exclude loops and parallel edges. If $\vert V\vert =\lambda$, we
call $X$ a {\sl $\lambda$-graph}, and whenever possible, we outright
assume that $V=\lambda$. A graph $X=(V,E)$ is {\sl $K(\kappa)$-free}, if
there is no clique of cardinal $\kappa$, i.e., $[T]^2\not\subseteq E$ holds
for every $T\in [V]^{\kappa}$. A {\sl $(\lambda,\kappa)$-graph} is a
$K(\kappa)$-free $\lambda$-graph. If $X_i=(V_i,E_i)$ ($i<2$) are graphs,
the one-to-one function $f\colon V_0 \to V_1$ is a {\sl weak (strong)
embedding} if $\{x,y\}\in E_0$ implies $\{f(x),f(y)\}\in E_1$ (if
$\{x,y\}\in E_0$ iff $\{f(x),f(y)\}\in E_1$). A {\sl weakly (strongly)
($\lambda,\kappa$)-universal graph} is a ($\lambda,\kappa$)-graph $X$ that
weakly (strongly) embeds every ($\lambda,\kappa$)-graph.
 
\medbreak
\centerline{\bf 1. When GCH holds}
 
\medskip
\proclaim Lemma 1. If $\lambda$ is strong limit,
$\lambda>\kappa\geq\omega$,  $\cf(\kappa)>\cf (\lambda)$
then there exists a strongly ($\lambda,\kappa$)-universal graph.

\medskip
\noindent{\bf Proof.} Let $\lambda=\sup\{\lambda_{\alpha}\colon \alpha<\cf
(\lambda)\}$, where the sequence is continuous, and
$2^{\lambda_{\alpha}}\leq \lambda_{\alpha+1}$, $\lambda_0=0$. Let $T$ be a
tree of height $\cf (\lambda)$ in which every $\alpha$-branch has
$\lambda_{\alpha+2}$ extensions on the $\alpha$-th level. Clearly, $\vert
T \vert =\lambda^{<\cf(\lambda)}=\lambda$. The vertex set of the universal graph
$X$ will be the disjoint union of some sets $\{A(t)\colon t\in T\}$ with
$\vert A(t)\vert =\lambda_{\alpha+1}$. No edge of $X$ will go between
$A(t)$ and $A(t')$ when $t$, $t'$ are incomparable in $T$. By induction on
$\alpha< \cf(\lambda)$, we determine for each $t\in T$ of height $\alpha$
how to build $X$ on $A(t)$, and how to join the vertices of $A(t)$ into 
$\bigcup \{A(t')\colon t'<t\}$. This latter set is of cardinal
$\lambda_{\alpha}$, with a graph on it, and we make sure that it will be
extended to a set of cardinal $\lambda_{\alpha+1}$, i.e., to some $A(t)$, in
all possible ways, such that the graph on $A(t)$ is $K(\kappa)$-free. This
is possible, as for every branch we have enough extensions reserved. It is
immediately seen that every ($\lambda,\kappa$)-graph embeds into $X$, one
only has to select the right branch.
 
The vertex set is of cardinal $\leq \vert T\vert \lambda=\lambda$. Finally,
a $K(\kappa)$ could only be produced along  a branch $\{A(t)\colon t\in
b\}$, but as $\vert b\vert \leq \cf(\lambda)<\cf(\kappa)$, some $A(t)$ must
contain a $K(\kappa)$, a contradiction, i.e., $X$ is a
($\lambda,\kappa$)-graph.
\medbreak
\proclaim Lemma 2. (F. Galvin) If $\lambda^{<\kappa}=\lambda$, then there is
no weakly ($\lambda,\kappa$)-universal graph.

\medskip
\noindent{\bf Proof.} Assume that $X=(\lambda,E)$ is
($\lambda,\kappa$)-universal. Let $Y=(V,G)$ be the following graph. The
elements of $V$ are those functions $f$ with $\Dom (f) <\kappa$ such that
$\Ran (f)$ is a clique in $E$. $\{f,g\}\in G$ iff $f \subset g$. Clearly,
$\vert V\vert =\lambda^{<\kappa}=\lambda$. If $\{f_{\alpha}\colon
\alpha<\kappa\}$ form a $K(\kappa)$, then they are compatible functions,
and their union $f=\bigcup\{f_{\alpha}\colon \alpha<\kappa\}$ injects
$\kappa$ into a clique of $X$, a contradiction, as $X$ is
$K(\kappa)$-free.
 
Assume that $g\colon V\to \lambda$ is a weak embedding of $Y$ into $X$. By
induction on $\alpha<\kappa$ we define $x_{\alpha}<\lambda$, $f_{\alpha}\in
V$ such that for $\beta<\alpha$ $\{x_{\beta},x_{\alpha}\}\in E$,
$f_{\beta}\subset f_{\alpha}$ (so $\{f_{\beta},f_{\alpha}\}\in G$) should
hold. If we succeed, we are done, as $\{x_{\alpha}\colon \alpha<\kappa\}$
is a clique again. If $\{x_{\beta},f_{\beta}\colon \beta<\alpha\}$ are
defined, let $f_{\alpha}$ be the following function: $\Dom
(f_{\alpha})=\alpha$, $f_{\alpha}(\beta)=x_{\beta}$ ($\beta<\alpha$).
$f_{\alpha}\in V$, as its range, $\{x_{\beta}\colon \beta<\alpha\}$ is a
clique. Put $x_{\alpha}=g(f_{\alpha})$. As by the way $f_{\alpha}$ is
constructed, $f_{\beta}\subset f_{\alpha}$ ($\beta<\alpha$), and $g$ is a
weak embedding, $x_{\alpha}$ will indeed, be joined into $x_{\beta}$ for
$\beta < \alpha$, and so the inductive step is successfully completed.
\medskip
\noindent{\bf Lemma 3.}{\sl\ If $\lambda$ is strong limit, $\kappa\leq
\lambda$, $\cf (\kappa)\leq \cf (\lambda)$, then there is no weakly
($\lambda,\kappa$)-universal graph.}
\medskip
\noindent{\bf Proof.} We can assume that $\kappa>\cf (\lambda)$, as
otherwise Lemma 2 gives the result. Assume that $X=(\lambda,E)$ is
($\lambda,\kappa$)-universal. Let $\{\kappa_{\alpha}\colon \alpha<\cf
(\kappa)\}$ be an increasing sequence of regular cardinals, cofinal in $\kappa$,
with $\kappa_0> \cf (\lambda)$. Let $F$ be the set of those $f$ functions
which satisfy the following requirements. $\Dom (f) < \cf (\kappa)$, for
$\alpha\in \Dom (f)$, $f(\alpha)$ is a bounded subset of $\lambda$ with
$\vert f(\alpha)\vert =\kappa_{\alpha}$, and $\bigcup \{ f(\alpha)\colon
\alpha< \Dom(f)\}$ is a clique in $X$. Let $V$, the vertex set of the
graph $Y=(V,G)$ be the disjoint union of the sets $\{ A(f)\colon f\in F\}$
where $\vert A(f)\vert =\kappa_{\alpha}$ if $\Dom(f)=\alpha$. Two distinct
vertces are joined iff one of them is in $A(f)$ the other in $A(f')$ for
some $f\subseteq f'$.

 Clearly, $\vert V\vert \leq \kappa \vert F \vert
=\lambda$. Assume that $T$ spans a clique in $Y$ and $\vert T\vert
=\kappa$. Then $T\subseteq \bigcup \{ A(f_{\gamma})\colon \gamma\in
\Gamma\}$ for a collection of pairwise compatible $f_{\gamma}$'s. $\sup
(\Dom(f_{\gamma}))=\cf(\kappa)$ as otherwise $\vert T\vert <\kappa$, but
then $\bigcup\{\Ran(f_{\gamma})\colon \gamma\in \Gamma\}$ is a $K(\kappa)$
in $X$, a contradiction. We therefore established that $Y$ is a
($\lambda,\kappa$)-graph.
 
Assume that $g\colon V\to \lambda$ is a weak embedding of $Y$ into $X$. By
induction on $\alpha <\cf (\kappa)$ we are going to define $f_{\alpha}\in
F$ such that $\Dom(f_{\alpha})=\alpha$, $f_{\alpha+1}(\alpha)\subseteq
g''A(f_{\alpha})$, and $f_{\beta}\subset f_{\alpha}$ whenever
$\beta<\alpha$. If this can be carried out, we reached a contradiction as
then $\bigcup \{\Ran(f_{\alpha})\colon \alpha< \cf(\kappa)\}$ is a
$K(\kappa)$ in $X$. There is no problem with the definition of $f_{\alpha}$
if $\alpha=0$ or limit. Assume that $f_{\alpha}$ is given.
$g''A(f_{\alpha})$ is a clique in $X$ of size $\kappa_{\alpha}=\cf(\kappa_{\alpha})>\cf(\lambda)$, so, there is a bounded (in $\lambda$) subset of it of
cardinal $\kappa_{\alpha}$, say, $S$. We can now define
$f_{\alpha+1}(\alpha)=S$, $f_{\alpha+1}(\beta)=f_{\alpha}(\beta)$
($\beta<\alpha$), the vertices in $f_{\alpha}(\beta)$ will be joined to
$S$, as by condition, $f_{\alpha}(\beta)=f_{\beta+1}(\beta)\subseteq g''
A(f_{\beta})$, $A(f_{\beta})$ is joined to $A(f_{\alpha})$ by the condition
$f_{\beta}\subset f_{\alpha}$, and $g$ is a weak embedding.
 
\medskip
From the known results and Lemmas 1--3 we can deduce the following.
\medskip
 
\noindent{\bf Theorem 1.}{\sl\ (GCH) Given $\lambda\geq \kappa$,
$\lambda\geq \omega$, there is a weakly/strongly
($\lambda,\kappa$)-universal graph iff $\kappa<\omega$ or $\cf (\kappa)>
\cf (\lambda)$.}
\medbreak
\centerline{\bf 2. The structure of the class of
$(\lambda,\kappa)$-graphs}
\medskip
\noindent In this Section we investigate the complexity of the class of
($\lambda,\kappa$)-graphs when there is no universal element in it. 
\medskip
\noindent{\bf Definition.} For $\lambda\geq \kappa$, ${\rm
CF}(\lambda,\kappa)$ is the minimal cardinal $\mu$ such that there is a family
$\{X_{\alpha}\colon \alpha<\mu \}$ of $(\lambda,\kappa)$-graphs, with the
property that every ($\lambda,\kappa$)-graph is weakly embedded into some
$X_{\alpha}$. ${\rm CF}^+(\lambda,\kappa)$ is the same with strong
embeddings. 
\medskip
 
Clearly, ${\rm CF}(\lambda,\kappa)\leq {\rm CF}^+(\lambda,\kappa)\leq
2^{\lambda}$. Also, ${\rm CF}(\lambda,\kappa)\leq \lambda$ iff ${\rm
CF}(\lambda,\kappa)=1$ iff there is a weakly ($\lambda,\kappa$)--universal
graph, and likewise for ${\rm CF}^+(\lambda,\kappa)$.
 
It was observed in [1] that ${\rm CF}^+(\omega,\omega)=\omega_1$. We
slightly extend that result.
\medbreak
\noindent{\bf Theorem 2.}{\sl\ If $\lambda\geq \kappa$, $\lambda$ is either
strong limit or of the form $\lambda=\mu^+=2^{\mu}$, $\cf(\kappa)=\omega$,
then ${\rm CF}^+(\lambda,\kappa)=\lambda^+$.}
\medskip
 
\noindent{\bf Proof.} From Lemmas 2--3, ${\rm CF}(\lambda,\kappa)\geq
\lambda^+$. Fix an increasing sequence $\kappa_n\to \kappa$, $\kappa_0=0$.
Call a structure $(A,<,X,R)$ a {\sl ranked graph} if $(A,<)$ is a
well--ordered set, $X$ is a graph on $A$, and $R$ is a function mapping
those bounded cliques of $X$ with order--type some $\kappa_n$ into the
ordinals, with the property that if clique $C'$ end--extends clique $C$, then
$R(C')<R(C)$. Obviously, then $X$ will be $K(\kappa)$--free. On the other
hand, if a $K(\kappa)$--free graph $X$ is given on a well--ordered set
$(A,<)$, then the tree 
$$T(X)=\{C\subseteq A\colon {\rm type}(C)=\kappa_n\ ({\rm some\ }n),
C{\rm\ clique\ }\}$$
endowed with end-extension, as the partial order, will be
$\omega$-branchless, so an ordinal valued function $R$ as above exists. If
$\vert A\vert =\lambda$, then $\vert T\vert = \lambda$, so only $\lambda$
ordinals are used, therefore $R(0)<\lambda^+$ holds. We call the minimal
possible $R(0)$ the rank of $X$.
 
Assume first that $\lambda$ is strong limit. Fix a continuous, cofinal
sequence $\{\lambda_{\alpha}\colon \alpha<\cf(\lambda)\}$ of cardinals with
$\lambda_0=0$ and $2^{\lambda_{\alpha}}\leq \lambda_{\alpha+1}$. 
 
For every $\xi <\lambda^+$ we are going to construct a graph that embeds
all graphs with rank $\xi$. 
 
Let $T$ be a tree with height $\cf(\lambda)$, with one root, such that whenever
$\alpha<\cf(\lambda)$, then every $\alpha$-branch has $\lambda_{\alpha+2}$
extensions to the $\alpha$-th level. For $t\in T$ on the $\alpha$-th
level, let $A(t)$ be an ordered set of order-type $\lambda_{\alpha+1}$,
such that the sets $\{A(t)\colon t\in T\}$ are pairwise disjoint. The
vertex set $V$ of our graph will be the union $V$ of these sets. We
partially order $V$ by assuming $A(t)<A(t')$ for $t<t'$, i.e., all elements
of $A(t)$ precede all elements of $A(t')$. 
 
For every $t\in T$, put $B(t)=\bigcup\{A(t')\colon t'<t\}$. By induction on
the height of $t$ we define $S(t)$, a ranked graph with ranks $\leq \xi$ on
$B(t)\cup A(t)$ such that if $b$ is an $\alpha$-branch, then all possible
end-extensions (if there are any) of the already defined structure on
$\bigcup\{A(t)\colon t\in b\}$ actually occur. This is possible, as there
are enough extensions of $b$ to the $\alpha$-th level. 
 
It is now obvious that all $(\lambda,\kappa)$-graphs of rank $\leq\xi$
embed into our tree. One only has to select the appropriate branch through
$T$. Also, $\vert V\vert =\vert T\vert
\lambda=\lambda^{<\cf(\lambda)}=\lambda$.We need to show that there is no $K(\kappa)$ in the
resulting graph. Assume that $U$ is a clique, $\vert U\vert =\kappa$. As we
joined vertices only in comparable $A(t)$'s, $U\subseteq
\bigcup\{A(t)\colon t\in b\}$ for some branch $b$. For some $t_n\in b$
($n=0,1,\dots$), it is true that the first $\kappa_n$ elements of $U$ are
bounded in $S(t_n)$, so they get a decreasing sequence of ordinals as
ranks, a contradiction.
 
The case $\lambda=\mu^+=2^{\mu}$ is actually simpler, we need one--element
$A(t)$'s, and having $\mu^+$ extensions of every branch of length $<\mu^+$.
 
\medbreak
 
Finally we show that under $\kappa^{<\kappa}=\kappa$, ${\rm
CF}(\kappa,\omega_1)$ can be as small as $\kappa^+$, and as large as
$2^{\kappa}$, and this latter value as large as we wish.
\def\AA{{\cal A}}
\medskip
\noindent{\bf Theorem 3.}{\sl\ Assume that in $V$, a model of GCH, $\mu$,
$\kappa>\omega$ are cardinals, $\cf(\mu)>\kappa=\cf(\kappa)$, then in a
cardinal and cofinality preserving forcing extension $V^P$, the GCH holds
below $\kappa$ and ${\rm CF}(\kappa,\omega_1)=2^{\kappa}=\mu$.}
\medskip
\noindent{\bf Proof.} If $\kappa=\lambda^+$, with $\cf (\lambda)=\omega$,
then we first add a $\square_{\lambda}$-sequence, i.e., a sequence
$\{C_{\alpha}\colon \alpha<\kappa,{\rm\ limit}\}$ with the following
properties:
\item {(1)} $C_{\alpha}\subseteq \alpha$ is closed, unbounded ;
\item {(2)} if $\gamma$ is a limit point of $C_{\alpha}$, then $C_{\gamma}=
\gamma\cap C_{\alpha}$ ;
\item {(3)} $\vert C_{\alpha}\vert < \kappa$.
 
It is well known that such a sequence can be added by a cardinal and
cofinality preserving forcing of size $\kappa$, so we may assume that it
exists in $V$. Fix such a sequence, and a sequence of cardinals
$\lambda_n\to \lambda$, and a one-to-one mapping $\phi_{\alpha,\beta}\colon
[\alpha,\beta)\to \lambda$ for each $\alpha<\beta<\kappa$. 
 
We call a countable set $A\subseteq \kappa$  $\sl low$, if ${\rm tp}(A)$ is
limit, and, if we put
 $\delta=\sup(A)$, $C_{\delta}=\{c_{\xi}\colon \xi<{\rm
tp}(C_{\delta})\}$ the increasing enumeration of $C_{\delta}$, then for
some $n<\omega$, $\phi_{c_{\xi},c_{\xi+1}}(a)<\lambda_n$ holds for $a\in
A$, $c_{\xi}\leq a<c_{\xi+1}$. 
 
If $\kappa>\omega_1$ is not of the form $\kappa=\lambda^+$, with
$\cf(\lambda)=\omega$, then we call every countable subset of limit type
low. 
\medskip
 
\noindent{\bf Claim 1.}{\sl\ The number of low subsets of some
$\alpha<\kappa$ is $<\kappa$.}
\medskip
\noindent{\bf Proof.} If $\kappa$ is not of the form $\lambda^+$ with
$\cf(\lambda)=\omega$, then $\vert \alpha\vert ^{\omega}<\kappa$. In the
other case the statement follows from property (3).
\medbreak
\noindent{\bf Claim 2.}{\sl\ If $B\subseteq \kappa$ is of order-type $\omega_1$, then for
some cofinal subset $B'\subseteq B$ it is true that if $\gamma<\sup(B')$ is
a limit point of $B'$, then $B'\cap \gamma$ is low.}
\medskip
\noindent{\bf Proof.} Put $\delta=\sup(B)$. Shrink $B$ to a cofinal
$B'\subseteq B$, such that the elements of $B'$ are separated by
$C_{\delta}$, and there is an $n<\omega$, such that if $c_{\xi}\leq b
<c_{\xi+1}$ for some $\xi$, then $\phi_{c_{\xi},c_{\xi+1}}(b)<\lambda_n$
($b\in B'$). Then the Claim follows from property (2) of the
$\square$-sequence.
 
If $\kappa$ is not of the form $\kappa=\lambda^+$ with
$\cf(\lambda)=\omega$ the choice $B'=B$ works.
 
\medbreak
The poset $(P,\leq)$ of the proof of the Theorem will be the $<\kappa$ support product
of $\mu$ copies of some poset $(Q,\leq)$ to be described below.
\medskip
$q\in Q$ if $q=(\delta,X,\AA)$ where $\delta<\kappa$, $X\subseteq
[\delta]^2$, $X$ is $K(\omega_1)$-free, if $\kappa>\omega_1$ the $\AA$ is a
family of low subsets of $\delta$, if $\kappa=\omega_1$, then $\AA$ is a
countable family of countable subsets of $\delta$ of limit type.
 Moreover, we require that if
$A\in\AA$, $\sup(A)\leq x<\delta$, then $A\times\{x\}\not\subseteq X$. 
 
$q'=(\delta',X',\AA')\leq q=(\delta,X,\AA)$ iff $\delta'\geq \delta$,
$X=X'\cap [\delta]^2$, $\AA=\AA'\cap [\delta]^{\aleph_0}$.
 
\medbreak
\noindent{\bf Claim 3.}{\sl\ $\vert Q\vert = \kappa$.}
\medskip
\noindent{\bf Proof.} For every $\delta<\kappa$ there are at most $\kappa$
many possibilities of selecting $X$, $\AA$ such that $(\delta,X,\AA)\in
Q$.
\medbreak
\noindent{\bf Claim 4.}{\sl\ Forcing with $(Q,\leq)$ does not introduce new
sequences of ordinals of length $<\kappa$.}
\medskip
\noindent{\bf Proof.} If $\kappa=\omega_1$, then $(Q,\leq)$ is
$<\omega_1$-closed. 

If $\kappa>\omega_1$, assume that $q\forc f\colon \tau
\to {\rm OR}$, $\tau<\kappa$. We construct the decreasing sequence of
conditions $\{q_{\alpha}=(\delta_{\alpha},X_{\alpha},\AA_{\alpha})\colon \alpha\leq \tau\}$ such that $q_0=q$, $q_{\alpha+1}\forc f(\alpha)=g(\alpha)$, and
if $\alpha$ is limit, then $\delta_{\alpha}=\sup\{\delta_{\beta}\colon
\beta<\alpha\}$, $X_{\alpha}=\bigcup\{X_{\beta}\colon \beta<\alpha\}$. If
$\cf(\alpha)\neq \omega$ then $\AA_{\alpha}=\bigcup\{\AA_{\beta}\colon
\beta<\alpha\}$, otherwise we add all the low subsets that are cofinal in
$\delta_{\alpha}$, 
to $\AA_{\alpha}$, as well. If we can carry out the construction, we are
done, $q_{\tau}$ determines all values of $f$. The only problem is if some
of the $X_{\alpha}$'s is not $K(\omega_1)$-free. Let $\alpha\leq \tau$ be
minimal such that there exists an uncountable 
 clique  $T\subseteq \delta_{\alpha}$. Clearly, $\cf(\alpha)=
\omega_1$. For some cofinal $T'\subseteq T$, if $\gamma< \delta_{\alpha}$
is a limit point of $T'$, then $T'\cap \gamma$ is low. There is a limit
$\beta<\alpha$ such that $\delta_{\beta}$ is a limit point of $T'$, so by
our construction $T'\cap \delta_{\beta}\in \AA_{\beta}$, so $T'\cap
\delta_{\beta}$ may not have been later extended to an $\omega_1$-clique.
\medbreak
\noindent{\bf Claim 5.}{\sl\ Forcing with $(P,\leq)$ does not introduce new
sequences of ordinals of length $<\kappa$.}
\medskip
\noindent{\bf Proof.} Similar to the previous proof.
\medbreak
\noindent{\bf Claim 6.}{\sl\ $(P,\leq)$ is $\kappa^+$-c.c.}
\medskip
\noindent{\bf Proof.} By Claim 3 and $\Delta$-system arguments.
\medbreak
If, in $V^P$, ${\rm CF}(\kappa,\omega_1)<\mu$, then a family of graphs
witnessing this is in a $<\mu$ sized subproduct of $P$. By the product
lemma we only need to show that forcing with $(Q,\leq)$ introduces a
$(\kappa,\omega_1)$-graph that cannot be embedded into any ground model
$(\kappa,\omega_1)$-graph. If $G\subseteq Q$ is generic, put
$Y=\bigcup\{X\colon (\delta,X,\AA)\in G\}$. 
\medbreak
\noindent{\bf Claim 7.}{\sl\ $Y$ is $K(\omega_1)$-free.}
\medskip
\noindent{\bf Proof.} If $\kappa=\omega_1$, $q\forc T$ is an
$\omega_1$-clique, select a decreasing sequence $q=q_0\geq q_1\geq\dots$
such that $q_{n+1}=(\delta_{n+1},X_{n+1},\AA_{n+1})\forc t_n\in T$,
$\delta_n < t_n <\delta_{n+1}$, and then put $q'=(\delta,X
,\AA)$ where $\delta=\lim \delta_n$, $X=\bigcup\{X_n:n<\omega\}$, and
$\AA=\bigcup\{\AA_n\colon n<\omega\}\cup \{\{t_n\colon
n<\omega\}\}$. Then $q' \forc T\subseteq  \delta$, a contradiction.
 
If $\kappa>\omega_1$, then by Claim 4 some $q=(\delta,X,\AA)$ determines all
elements of $T$, the alleged $\omega_1$-clique. We can assume that
$T\subseteq \delta$, but then $X$ is not $K(\omega_1)$-free, a
contradiction.
\medskip
\noindent{\bf Claim 8.}{\sl\ $Y$ does not embed into any ground model
$(\kappa,\omega_1)$-graph.}
\medskip
\noindent{\bf Proof.}
Assume that $q\forc f:\kappa\to\kappa$ is an embedding of $Y$ into some
ground model $(\kappa,\omega_1)$-graph, $Z$. By induction on
$\alpha<\omega_1$ construct the decreasing sequence
$q_{\alpha}=(\delta_{\alpha},X_{\alpha},\AA_{\alpha})$ such that $q_0=q$,
$q_{\alpha+1}\forc f(\delta_{\alpha})=g(\alpha)$, for $\alpha$ limit
$\delta_{\alpha}=\lim \{\delta_{\beta}:\beta<\alpha\}$,
$X_{\alpha}=\bigcup\{X_{\beta}:\beta<\alpha\}$,
$\{\delta_{\beta},\delta_{\alpha}\}\in X_{\alpha+1}$ for $\beta<\alpha$,
and $\AA_{\alpha}=\bigcup\{\AA_{\beta}:\beta<\alpha\}$. The only problem
with the definition would be that $A\subseteq
\{\delta_{\beta}:\beta<\alpha\}$ for some $A\in \AA_{\alpha}$. But then,
$\sup(A)$ is of the form $\delta_{\gamma}$ for some limit $\gamma\leq
\alpha$, and no set of that form was added to $\AA_{\gamma}$. 
 
We can therefore define the sequence, but then the range of $g$ will be a
$K(\omega_1)$ in Z, a contradiction.
 
\def\ZZ{{\cal Z}}
\medskip
\noindent{\bf Theorem 4.}{\sl\ If, in a model of GCH, $\mu$,
$\kappa>\omega$ are cardinals, with $\cf(\mu)>\kappa=\cf(\kappa)$, then, in
some cardinal and cofinality preserving extension the GCH holds below
$\kappa$, $2^{\kappa}=\mu$, and ${\rm CF}^+(\kappa,\omega_1)=\kappa^+$.}
\medskip
\noindent{\bf Proof.} Again, as in the proof of Theorem 3, we can assume,
that if $\kappa=\lambda^+$, with $\lambda>\cf(\lambda)=\omega$, then
$\square_{\lambda}$ holds in the ground model. We also assume that the GCH
holds below $\kappa$ and $2^{\kappa}=\mu$.
 
In a $<\kappa$-support iteration of length $\kappa^+$, we add a family
witnessing ${\rm CF}^+(\kappa,\omega_1)=\kappa^+$. Factor $Q_{\alpha}$ will
add a $(\kappa,\omega_1)$-graph that strongly embeds every
$(\kappa,\omega_1)$-graph of $V^{P_{\alpha}}$. Notice, that if the forcing
does not collapse cardinals, then $\square_{\lambda}$ will still hold at
every stage.

We first define and investigate one step of the iteration.
 
Let $(Q,\leq)$ be the following poset. $q=(\delta, X, \AA,\ZZ,F)\in Q$, if
$\delta<\kappa$, $X\subseteq [\delta]^2$ is a $K(\omega_1)$-free graph,
$\AA\subseteq [\delta]^{\aleph_0}$ is a family of low sets
($\kappa>\omega_1$), is a countable family of limit type subsets of
$\delta$ ($\kappa=\omega_1$). $\ZZ$ is a family of $<\kappa$ many
$(\kappa,\omega_1)$-graphs, $F:\ZZ\times \delta \to \delta$ is a function
such that if $Z\in \ZZ$ then the mapping $x\mapsto F(Z,x)$ is a strong
embedding of $Z\vert \delta$ into $X$, and the following two more
conditions hold.
 
\item {(1)} If $A\in\AA$, $\sup(A)\leq x<\delta$, then
$A\times\{x\}\not\subseteq X$ ;
\item {(2)} if $A\in\AA$, $Z\in\ZZ$, then
$A\not\subseteq F''(\{Z\}\times\delta)$.

$q'=(\delta',X',\AA',\ZZ',F')\leq q=(\delta,X,\AA,\ZZ,F)$ if
$\delta'\geq\delta$, $X=X'\cap[\delta]^2$, $\ZZ'\supseteq \ZZ$,
$\AA=\AA'\cap [\delta]^{\aleph_0}$ and, moreover, 
 
\item {(3)} if $Z_0\neq Z_1\in \ZZ$, $\delta\leq x,y<\delta'$, then
$F'(Z_0,x)\neq F'(Z_1,y)$.
 
\medbreak
\noindent{\bf Claim 1.}{\sl\ $(Q,\leq)$ is transitive.}
\medbreak
\noindent{\bf Proof.} Assume that $q_0\geq q_1\geq q_2$,
$q_i=(\delta_i,X_i, \AA_i,\ZZ_i,F_i)$ ($i<3$). In establishing $q_0\geq
q_2$ only condition (3) could cause problems, but it will not: if $Z_0\neq
Z_1\in \ZZ_0$, $\delta_0\leq x<\delta_1\leq y<\delta_2$, then
$F_2(Z_0,x)\neq F_2(Z_1,y)$ as the first element is in
$[\delta_0,\delta_1)$, the second is in $[\delta_1,\delta_2)$.
 
\medbreak
\noindent{\bf Claim 2.}{\sl\ If $\varepsilon<\kappa$,
$D=\{(\delta,X,\AA,\ZZ,F):\delta\geq \varepsilon\}$ is dense.}
\medskip
\noindent{\bf Proof.} We can extend a given $(\delta,X,\AA,\ZZ,F)$ to a
large enough $\delta'$ by mapping $Z\vert [\delta,\delta')$ ($Z\in\ZZ$) onto disjoint
sets, not extending $\AA$, $\ZZ$, and adjusting $X$. Conditon (1) won't
cause problem, as by (2) no $A\in\AA$ will be forced to be joined to a
vertex.
 
\medskip
\noindent{\bf Claim 3.}{\sl\ If $Z$ is a $(\kappa,\omega_1)$-graph, then
$D=\{(\delta,X,\AA,\ZZ,F):Z\in \ZZ\}$ is dense.}
\medskip
\noindent{\bf Proof.} A similar argument works.
\medbreak
\noindent{\bf Claim 4.}{\sl\ Forcing with $(Q,\leq)$ doesn't introduce sequences
 of ordinals of length
$<\kappa$.}
\medskip
\noindent{\bf Proof.} $(Q,\leq)$ is $<\omega_1$-closed, and this is enough
if $\kappa=\omega_1$.
 
Assume that $\kappa>\omega_1$. Let $q\forc f:\tau\to {\rm OR}$, $\tau
<\kappa$. By induction on $\alpha\leq \tau$ we define the decreasing
sequence $\{q_{\alpha}=(\delta_{\alpha},X_{\alpha},\AA_{\alpha},\ZZ_{\alpha},F_{\alpha}):\alpha\leq\tau\}$ such that $q_{\alpha+1}\forc
f(\alpha)=g(\alpha)$, and for limit $\alpha$,
$\delta_{\alpha}=\sup\{\delta_{\beta}:\beta<\alpha\}$, 
$X_{\alpha}=\bigcup\{X_{\beta}:\beta<\alpha\}$, 
$\ZZ_{\alpha}=\bigcup\{\ZZ_{\beta}:\beta<\alpha\}$,
 $F_{\alpha}=\bigcup\{F_{\beta}:\beta<\alpha\}$. If $\cf(\alpha)>\omega$,
we take $\AA_{\alpha}=\bigcup\{\AA_{\beta}:\beta<\alpha\}$, otherwise we
add all cofinal in $\delta_{\alpha}$ low subsets $A$, for wich there is no
$Z\in \ZZ_{\alpha}$ with $A\subseteq F_{\alpha}''(\{Z\}\times\delta_{\alpha})$.
The only thing we have to show is that no $K(\omega_1)$ will be created. We
may assume, that $\alpha\leq \tau$ is limit, $T\subseteq \delta_{\alpha}$
is cofinal, and $T$ is an uncountable clique in $X_{\alpha}$. We can assume that
segments of $T$ of limit type are low sets. 
 As $T$ could grow, for a
club subset $C\subseteq \alpha$, of order type $\omega_1$, it is true that
if $\beta\in C$, then $T\cap\delta_{\beta}\subseteq
F_{\beta}''(\{Z\}\times\delta_{\beta})$ for some $Z\in \ZZ_{\beta}$. By
conditon (3), there can be only one such $Z$. If, moreover $\beta$ is a
limit point of limit points of $C$, then there is a $h(\beta)<\beta$, such
that for $h(\beta)<\gamma\leq \beta$ this $Z$ for $\gamma$ 
 is the same. By the pressing
down lemma, $h$ is bounded on an unbounded subset, so $T\cap
\delta_{\beta}\subseteq F_{\alpha}''(\{Z\}\times \delta_{\beta})$ for
uncountably many $\beta<\alpha$, but then the inverse image of $T$ will be
a $K(\omega_1)$ in $Z$, a contradiction.
 
\medbreak
Let $Y$ be the  graph added by $Q$, i.e., if $G\subseteq Q$ is
generic, then $Y=\bigcup\{X:(\delta,X,\AA,\ZZ,F)\in G\}$.
\medbreak
\noindent{\bf Claim 5.}{\sl\ $Y$ is $K(\omega_1)$-free.}
\medskip
\noindent{\bf Proof.} If $\kappa=\omega_1$, $q\forc T$ is an
$\omega_1$-clique in $Y$, then an argument as above shows that there is a
decreasing sequence $\{q_{\alpha}:\alpha<\omega_1\}$ determining more and
more elements of $T$, and we can freeze $T$ unless it is covered by
$\bigcup\{F_{\alpha}''(\{Z\}\times\delta_{\alpha}):\alpha<\omega_1\}$ for
some $Z$, which again gives a $K(\omega_1)$ in $Z$.
 
If $\kappa>\omega_1$, by the above Claim, the supposed clique $T$ is in the
ground model, some $q\in G$ contains in its $X$-part, a contradiction.
\medbreak
The iteration $(P_{\alpha},Q_{\alpha}:\alpha\leq \kappa^+)$ is defined as a
$<\kappa$-support iteration, with $Q_{\alpha}$ as the above $Q$, defined in
$V^{P_{\alpha}}$. 
 
\medbreak
In $Q_{\alpha}$, let $D_{\alpha}$ be the set of those conditons of the form
$q=(\delta,X,\AA,\ZZ,F)$ for which it is true that $Z_0\neq Z_1\in\ZZ$ implies
that $Z_0\vert \delta\neq Z_1\vert \delta$.
\medbreak
 
\noindent{\bf Claim 6.}{\sl\ $D_{\alpha}$ is dense in $Q_{\alpha}$.}
\medskip
\noindent{\bf Proof.} Using Claim 1, with $\varepsilon$ large enough.
\medbreak
If $q=(\delta,X,\AA,\ZZ,F)
\in Q_{\alpha}$ we put $\ell(q)=(\delta,X,\AA,\ZZ\vert \delta,F)$. Let
$E_{\alpha}$ be the following subset of $P_{\alpha}$. $p\in P_{\alpha}$ if
for all $\beta<\alpha$, $p\vert \beta$ determines $\ell(p(\beta))$ and
forces that $p(\beta)\in D_{\beta}$.
\medbreak
\noindent{\bf Claim 7.}{\sl\ For every $\alpha\leq \kappa^+$
\item {(a)}  $E_{\alpha}$ is dense in $P_{\alpha}$
;
\item {(b)} forcing with $P_{\alpha}$ does not add sequences of ordinals of
length $<\kappa$.}
\medskip
\noindent{\bf Proof.} Assume first that $\kappa>\omega_1$.
The proof is by induction on $\alpha\leq \kappa^+$. If (b) holds
for $\alpha$, then it holds for $\alpha+1$, by Claim 4. Assume that (a) and
(b) hold for $\alpha$, and $p\in P_{\alpha+1}$. We may assume that $p\vert
\alpha\forc p(\alpha)\in D_{\alpha}$. As (b) holds for $\alpha$, there is a
$q\leq p\vert \alpha$ which determines $p(\alpha)$. Extend $q$ to an $r\in
E_{\alpha}$, then take $r\cup p(\alpha)\in E_{\alpha+1}$.
 
Assume that $\alpha$ is limit, $p\in P{\alpha}$. In order to prove  (a)
for $\alpha$, we may assume that $\supp(p)$ is cofinal in $\alpha$, let
$\{\alpha_{\xi}:\xi<\tau\}$ converge to $\alpha$. We define
$\{p_{\xi}:\xi<\tau\}$, a decreasing sequence of conditions. $p_0=p$. $p_{\xi}
\vert \alpha_{\xi}\in E_{\alpha_{\xi}}$, and $p_{\xi}\leq p_{\zeta}$,
$p_{\xi}\vert [\alpha_{\xi},\alpha)=p_{\zeta}\vert [\alpha_{\zeta},\alpha)$
hold for $\zeta<\xi$. If $\xi$ is limit, $\beta\geq\alpha_{\xi}$, the names
$p_{\zeta}(\beta)$ are identical, so we can take it as $p_{\xi}(\beta)$. If
$\beta<\alpha_{\xi}$, we take $p_{\xi}(\beta)$ as
$\bigcup\{p_{\zeta}(\beta)\}$ by adding all low subsets which can be added,
as in Claim 4. We show that $p_{\xi}$ is a condition. To this end, we show
by induction on $\beta<\alpha$ that $p_{\xi}\vert \beta$ is a condition.
The limit case is trivial. The problem with $p_{\xi}(\beta)$ can only be
that its $X$ part contains a $K(\omega_1)$, but then, as in the proof of
Claim 4, we get that $p_{\xi}\vert \beta\forc Z$ is not $K(\omega_1)$-free
for some $Z\in \ZZ$. 
 
If $\alpha$ is limit and we are to show (b) for $\alpha$, and $p\forc f:\tau\to
{\rm OR}$ for some $\tau<\kappa$, we can define a decreasing, continuous
sequence $\{p_{\xi}:\xi \leq\tau\}$ with $p_{\xi}\forc f(\xi)=g(\xi)$,
$p_{\xi}\in E_{\alpha}$. This can be carried out, as above, and then
$p_{\tau}$ decides $f$.

For $\kappa=\omega_1$, (b) follows from the fact that we iterate a
countably closed poset with countable supports, and for (a) an easy
inductive proof can be given, as for the other case above.
 
\medbreak
\noindent{\bf Claim 8.}{\sl\ $P_{\kappa^+}$ is $\kappa^+$-c.c.}
\medskip
\noindent{\bf Proof.} Given $\kappa^+$ conditions, we can assume that they
are from $E_{\kappa^+}$. By the usual $\Delta$-system arguments we can find
two of them $p$ and $p'$ such that $\ell(p(\alpha))=\ell(p'(\alpha))$ holds
for every $\alpha\in \supp(p)\cap \supp(p')$. We show that $p\cup p'$ is a
condition (though not necessarily in $E_{\kappa^+}$).
 
To this end, we show that $(p\cup p')\vert \alpha\in P_{\alpha}$ by
induction on $\alpha$. All cases are trivial, except when $\alpha=\beta+1$,
$\beta\in\supp(p)\cap \supp(p')$. What we have to show is that the $F$ part
of $(p\cup p')(\beta)$ is well-defined, i.e., if $Z=Z'$ are from the $\ZZ$
part, then $F(Z,x)=F(Z',x)$ ($x<\delta$). But this will hold (or, more
precisely, will be forced to hold by $(p\cup p')\vert\beta$) as $F(Z,x)$ is determined
by $Z\vert \delta$ and by $x$, and it is determined the same way in $p$ and
$p'$.
 
\medskip
From the last Claim, every $(\kappa,\omega_1)$-graph appears in some
intermediate extension, and so it is embedded into the next graph,
$Y_{\alpha}$, by $Q_{\alpha}$. We still have to show that $Y_{\alpha}$
remains $K(\omega_1)$-free under the further extensions. This follows from
Claim 7(b) if $\kappa>\omega_1$, and from the following statement which is
a special case of a well-known lemma about forcing.
\medskip
\noindent{\bf Claim 9.}{\sl\ If, in $V$, $Y$ is a $K(\omega_1)$-free graph,
$P$ is an $<\omega_1$-closed frocing, then, in $V^P$, $Y$ is still
$K(\omega_1)$-free.}
\medskip
\noindent{\bf Proof.} If $p\forc T$ is an uncountable clique, select
$\{p_{\alpha}:\alpha<\omega_1\}$ fixing more and more elements of $T$,
$p_0=p$.
\medskip
\noindent{\bf Remark.} With the technique of Theorem 4 it is possible to
show that if $\mu\geq \nu> \kappa$, $\cf(\mu)>\kappa$, and $\nu,\kappa$ are
regular, then it is consistent that $2^{\kappa}=\mu$,  ${\rm
CF}(\kappa,\omega_1)=\nu$, and GCH holds below $\kappa$. Add a sequence $\{Y_{\alpha}:\alpha<\nu\}$ 
, rather than of length $\kappa^+$, as in Theorem 4. One only
has to observe that $Y_{\alpha}$ does not embed into any $K(\omega_1)$-free
graph in $V^{P_{\alpha}}$, this can be proved similarly to Claim 8 in Theorem
3.
\goodbreak
\centerline{\bf References}
\medskip
\item {[1]} A.~Hajnal, P.~Komj\'ath: Embedding graphs into colored graphs,
{\sl Trans. of the Amer. Math. Soc. \bf 307} (1988), 395--409.
\item{[2]} P. Komj\'ath and J\'anos Pach, Universal elements and the complexity
of certain classes of infinite graphs,
  {\sl Discrete Math. \bf 95 }(1991) 255--270. 
\item{[3]} P.~Komj\'ath, J.~Pach: The complexity of a class of
infinite graphs,
{\sl Combinatorica}, to appear.
\item {[4]} R.~Rado: Universal graphs and universal functions, {\sl Acta
Arith., \bf 9} (1964), 331--340.
\item {[5]} R.~Rado: Universal graphs, in: {\sl A Seminar in Graph Theory},
(eds. Harary, Beineke), Holt, Rinehart, and Winston Co., 1967.
 
\bye